\documentclass[11pt]{amsart}
\usepackage{amsmath,amssymb,amsthm}
\usepackage[latin1]{inputenc}
\usepackage{version,tabularx,multicol}
\usepackage{graphicx,float}

\newcommand{\reff}[1]{(\ref{#1})}

\theoremstyle{plain}
\newtheorem{theo}{Theorem}[section]

\newtheorem{cor}[theo]{Corollary}
\newtheorem{prop}[theo]{Proposition}

\newtheorem{lem}[theo]{Lemma}

\theoremstyle{remark}
\newtheorem{REM}[theo]{Remark}
\newenvironment{rem}
{\vspace{\parskip}\begin{REM}}
{\hspace*{\fill}$\lozenge$\end{REM}}
\newtheorem{EXA}[theo]{Example}
\newenvironment{exa}
{\vspace{\parskip}\begin{EXA}}
{\hspace*{\fill}$\triangle$\end{EXA}}

\newcommand{\ca}{{\mathcal A}}

\newcommand{\cf}{{\mathcal F}}
\newcommand{\cg}{{\mathcal G}}

\newcommand{\cn}{{\mathcal N}}

\newcommand{\cp}{{\mathcal P}}
\newcommand{\cq}{{\mathcal Q}}

\newcommand{\E}{{\mathbb E}}

\newcommand{\N}{{\mathbb N}}
\renewcommand{\P}{{\mathbb P}}

\newcommand{\cQ}{{\mathcal Q}}

\newcommand{\cJ}{{\mathcal J}}
\newcommand{\cI}{{\mathcal I}}
\newcommand{\R}{{\mathbb R}}

\newcommand{\Var}{{\rm Var}}
\newcommand{\Cov}{{\rm Cov}}

\newcommand{\ind}{{\bf 1}}

\newcommand{\Card}{{\rm Card}\;}

\newcommand{\inv}[1]{\mathop{\frac{1}{ #1}}\nolimits}
\newcommand{\expp}[1]{\mathop {\mathrm{e}^{ #1}}}

\begin{document}

\title[WRMC]{Does waste-recycling really improve
 the multi-proposal Metropolis-Hastings Monte Carlo algorithm?}

\date{\today}

\author{Jean-François Delmas}

\address{
Jean-Fran\c cois Delmas and Benjamin Jourdain,
CERMICS, \'Ecole des Ponts, Université Paris-Est, 6-8
av. Blaise Pascal, 
  Champs-sur-Marne, 77455 Marne La Vallée, France. Research supported by
the ANR program ADAP'MC.}

\email{delmas@cermics.enpc.fr, jourdain@cermics.enpc.fr}

\author{Benjamin Jourdain}

\begin{abstract}
  The waste-recycling  Monte Carlo
  (WR)  algorithm introduced  by  physicists is  a  modification of  the
  (multi-proposal) Metropolis-Hastings algorithm, which makes use of all
  the proposals in 
  the empirical mean, whereas the standard (multi-proposal) Metropolis-Hastings algorithm
  only uses  the accepted proposals.  In this paper, we extend the WR
  algorithm into a general control variate technique and exhibit the
  optimal choice of the control variate in terms of asymptotic
  variance. We
  also give an  example which shows that in contradiction to the
  intuition of physicists, the WR
  algorithm can have  an asymptotic variance larger than  the one of the
  Metropolis-Hastings algorithm.  However, in the particular case of the
  Metropolis-Hastings  algorithm called  Boltzmann  algorithm, we  prove
  that   the   WR   algorithm   is  asymptotically   better   than   the
  Metropolis-Hastings algorithm. This last property is also true for the
  multi-proposal Metropolis-Hastings algorithm. In this last framework,
  we consider a linear parametric generalization of WR, and we propose an
  estimator of the explicit optimal parameter using the proposals. 

\end{abstract}



\maketitle

\section{Introduction}

The  Metropolis-Hastings algorithm  is used  to compute  the expectation
$\langle  \pi,f\rangle$ of a  function $f$  under a  probability measure
$\pi$  difficult to  simulate.   It  relies on  the  construction by  an
appropriate   acceptation/rejection   procedure   of  a   Markov   chain
$(X_k,k\geq 0)$ with transition kernel $P$ such that $\pi$ is reversible
with respect to $P$ and  the quantity of interest $\langle \pi,f\rangle$
is estimated by  the empirical mean $I_n(f)=\inv{n}\sum_{k=1}^n f(X_k)$.
We shall recall the well-known properties of this estimation (consistency,
asymptotic  normality) in  what follows.  In particular  the  quality or
precision of the algorithm is measured through the asymptotic variance of
the estimator of $\langle \pi, f \rangle$.

The   waste-recycling  Monte   Carlo  (WR)   algorithm,   introduced  by
physicists,  is  a modification  of  the Metropolis-Hastings  algorithm,
which makes use of all the  proposals in the empirical mean, whereas the
standard Metropolis-Hastings algorithm only uses the accepted proposals.
To  our knowledge,  the WR  algorithm was  first introduced  in  1977 by
Ceperley, Chester and Kalos  in equation (35) p.3085 \cite{cck}. Without
any proof,  they claim that ``The  advantage of using this  form is that
some information about  unlikely moves appears in the  final answer, and
the variance  is lowered''. It  is commonly assumed among  the physicists
and supported by  most of the simulations that the  WR algorithm is more
efficient than the Metropolis-Hastings algorithm, that is the estimation
given by  the WR  algorithm is consistent  and has a  smaller asymptotic
variance.  An other  way to speed up the  Metropolis-Hastings algorithm could  be to  use multiple
proposals  at  each step  instead  of  only  one. According  to  Frenkel
\cite{f1:wrmc}, the waste recycling can be particularly useful for these
algorithms where many states are rejected.

Our aim  is to clarify  the presentation of  the WR algorithms  with one
proposal  and with  multiple proposals  and to  present a  first rigorous
study of those algorithms.  We will give in Section \ref{sec:did}
an introduction  to our results in the finite state  space case.   Our
main new results are stated in Theorem \ref{theo:princ2},
which is   a first  step  towards  the
comparison  of   the  asymptotic  variances.   We   shall  detail  their
consequences in the  didactic Section \ref{sec:did} for:
\begin{itemize}
   \item[-] the  WR algorithm  through Propositions
     \ref{prop:CV} (consistency of the estimation),   \ref{asvar}
     (asymptotic normality)  and \ref{prop:MB} (a first partial answer
     to the initial question: Does waste-recycling really improve the Metropolis-Hastings Monte Carlo algorithm?),
   \item[-] the multi-proposal WR algorithm  through Propositions
     \ref{prointmul} (consistency of the estimation and asymptotic
     normality) and \ref{prop:mult-MB} (a second partial answer
     to the inital question: Does waste-recycling really improve
 the Metropolis-Hastings Monte Carlo algorithm?).
\end{itemize}

The study of  the WR estimator in the form  $I_n(f)+J_n(f)$, for a given
functional $J$,  leads us  to rewrite the  WR algorithm as  a particular
case of a general control  variate problem by considering the estimators
$I_n(f)+J_n(\psi)$ where the function  $\psi$ is possibly different from
$f$.  In the multi-proposal  framework, the consistency (or convergence)
of this  general algorithm  and its asymptotic  normality are  stated in
Theorem \ref{theo:princ2} in Section  \ref{sec:Egene}.  We also give its
asymptotic variance and prove that the optimal choice of $\psi$ in terms
of asymptotic  variance is  the solution, $F$,  of the  Poisson equation
\reff{Poisson}.   This  choice  achieves  variance  reduction,  but  the
function  $F$ is  difficult to compute. It is possible to replace it by an 
approximation.  In some  sense, $f$ is such an approximation and
for  this particular  choice we  recover the  Waste  Recycling estimator
introduced by physicists.  In Section \ref{sp} which is dedicated to the
single proposal  case, we give  a simple counter-example  (see paragraph
\ref{sec:ce})  which shows  that the  WR algorithm  does not  in general
improve the Metropolis-Hastings algorithm : the WR algorithm can have an
asymptotic  variance  larger than  the  one  of the  Metropolis-Hastings
algorithm.  Since,   Ath\`enes  \cite{a}  has   also  observed  variance
augmentation in some numerical  computations of free energy. However, in
the  particular   case  of  the   Metropolis-Hastings  algorithm  called
Boltzmann  algorithm, we  prove in  Section \ref{sec:Bolt}  that  the
(multi-proposal) WR
algorithm   is  asymptotically   better  than   the  (multi-proposal)
Metropolis-Hastings 
algorithm. In this particular framework, we explicit the optimal value
$b_\star$ of $b$ for the parametric control variate $J_n(bf)$. This optimal
value can be estimated using the Makov chain $(X_k, 0\leq k\leq n)$.

$ $

{\bf Acknowledgments.} We warmly thank Manuel Ath\`enes (CEA Saclay) for
presenting the waste recycling Monte Carlo algorithm to us and Randal
Douc (CMAP \'Ecole Polytechnique) for numerous fruitful discussions. We
also thank the referees for their valuable comments. 

\section{Didactic version of the results}
\label{sec:did}

For simplicity, we
assume in the present section that $E$ is a finite set. Let $\langle \nu,h \rangle=\sum_{x\in
  E} \nu(x) h(x)$ denote the ``integration'' of a real function  defined
on $E$, $h=(h(x), x\in E)$, w.r.t. to a measure on $E$, $\nu=(\nu(x),
x\in E)$.  

Let $\pi$ be  a probability measure on $E$ such  that $\pi(x)>0$ for all
$x\in   E$   and   $f$   a   real  function   defined   on   $E$.    The
Metropolis-Hastings  algorithm  gives an  estimation  of $\langle  \pi,f
\rangle$  as the  a.s.  limit  of the  empirical mean  of  $f$, $\inv{n}
\sum_{k=1}^n f(X_k)$, as $n$ goes  to infinity, where $X=(X_n, n\geq 0)$
is a  Markov chain which is  reversible with respect  to the probability
measure $\pi$.

\subsection{The Metropolis-Hastings algorithm} 
\label{sec:mh}

The Markov chain $X=(X_n, n\in \N)$ of the Metropolis-Hastings algorithm
is built  in the  following way.  Let  $Q$ be an  irreducible transition
matrix over  $E$ such  that for  all
$x,y\in E$, if $Q(x,y)=0$ then  $Q(y,x)=0$. The transition matrix $Q$ is
called the selection matrix.

For $x,y\in E$ such that $Q(x,y)>0$, let $(\rho(x,y),\rho(y,x))\in
(0,1]^2$ be such that
\begin{equation}
   \rho(x,y)\pi(x)Q(x,y)=\rho(y,x)\pi(y)Q(y,x).\label{prerevers}
\end{equation}
The function $\rho$ is viewed as an acceptance probability. For example,
one gets such a function $\rho$ by setting
\begin{equation}
   \rho(x,y)=
\gamma\left(\frac{\pi(y)Q(y,x)}{\pi(x)Q(x,y)}\right),  \quad\text{for
  all} \quad  x,y\in E \quad\text{s.t.}\quad Q(x,y)>0,
\label{casga} 
\end{equation}
where  $\gamma$  is  a  function   with  values  in  $(0,1]$  such  that
$\gamma(u)=u\gamma(1/u)$.  Usually, one takes  $\gamma(u)=\min(1,u)$ for
the Metropolis  algorithm. The case $\gamma(u)=u/(1+u)$ is  known as the
Boltzmann algorithm or Barker algorithm.


Let $X_0$  be a  random variable taking  values in $E$  with probability
distribution $\nu_0$.  At step  $n$, $X_0, \ldots,  X_n$ are  given. The
proposal at  step $n+1$, $\tilde  X_{n+1}$, is distributed  according to
$Q(X_n,\cdot)$.   This    proposal   is   accepted    with   probability
$\rho(X_n,\tilde X_{n+1})$  and then $X_{n+1}=\tilde X_{n+1}$.  If it is
rejected, then we set $X_{n+1}=X_n$.

It is easy to check that $X=(X_n, n\geq 0)$ is a Markov chain with
transition matrix $P$ defined by 
\begin{equation}
   \label{eq:def-P}
\forall x,y \in E, \quad P(x,y)=\begin{cases}
   Q(x,y) \rho(x,y)&\text{if $x\neq y$,}\\
 1- \sum_{z\neq x} P(x,z)&\text{if $x= y$.}
\end{cases}
\end{equation}
Furthermore $X$ is reversible w.r.t. to the probability measure $\pi$:
$\pi(x) P(x,y)=\pi(y) P(y,x)$ for all $x,y\in E$. This property is
also called detailed balance. By summation over
$y\in E$, one deduces that $ \pi$ is an
invariant probability for $P$ (i.e. $\pi P=\pi$). The irreducibility of $Q$ implies that $P$ is
irreducible. Since the probability measure $\pi$ is invariant for $P$, we
deduce that $X$ is  positive recurrent  with (unique)  invariant
probability measure
$\pi$. In  particular, for any real  valued function $f$  defined on $E$, the ergodic theorem
(see e.g. \cite{mt:mcss}) implies the consistency of the estimation: 
\[
\lim_{n\rightarrow \infty } I_n(f)=\langle \pi, f \rangle \quad\text{a.s.},
\]
where
\begin{equation}
   \label{eq:def-In}
I_n(f)=\inv{n} \sum_{k=1}^n f(X_k).
\end{equation}
The asymptotic normality of the estimator $I_n(f)$ is given by the following
central limit theorem (see
\cite{d:rim} or \cite{mt:mcss})
\[
\sqrt{n} \left(I_n(f)- \langle \pi,f \rangle\right) 
\; \xrightarrow[n\rightarrow \infty ]{(d)} \;
\cn(0,\sigma(f)^2).
\]
Here $\cn(0,\sigma^2)$  denotes the  Gaussian distribution with  mean 0
and  variance  $\sigma^2$, the  convergence holds  in the  distribution
sense and \begin{equation}
   \label{eq:def-s2}
\sigma(f)^2
=\langle\pi,F^2\rangle-\langle\pi,(PF)^2\rangle.  
\end{equation}
where $F$ denotes the unique solution up to an additive constant of the
Poisson equation
\begin{equation}
 \label{Poisson} 
 F(x)-PF(x)=f(x) -
\langle \pi, f \rangle, \quad x\in E
\end{equation}
and $Ph(x)=\sum_{y\in E}P(x,y)h(y)$.
Improving the Metropolis-Hastings algorithm means  exhibiting  other estimators of $\langle \pi,f \rangle$ that  are still
consistent (i.e. estimators which converge a.s.   to
$\langle \pi,  f \rangle$)  but  with an asymptotic variance  smaller than
$\sigma(f)^2$.   


\subsection{WR algorithm}
The  classical estimation of  $\langle \pi,f  \rangle$ by  the empirical
mean $I_n(f)$  makes no  use of the  proposals $\tilde{X}_k$  which have
been  rejected.  For  a  long  time, physicists  have  claimed that  the
efficiency of the estimation can be improved by including these rejected
states in  the sampling  procedure.  They suggest  to use  the so-called
Waste-Recycling Monte Carlo (WR)  algorithm, which consists in replacing
$f(X_k)$  in  $I_n(f)$  by   a  weighted  average  of  $f(X_{k-1})$  and
$f(\tilde{X}_k)$. For the natural choice of weights corresponding to the
conditional expectation of $f(X_k)$  w.r.t. $(X_{k-1}, \tilde X_k)$, one
gets the following estimator of $\langle \pi,f\rangle$:
\begin{align*}
   I^{WR}_n(f)
&=\inv{n} \sum_{k=0}^{n-1} \E\left[f(X_{k+1})|X_k, \tilde
  X_{k+1}\right]\\ 
&=\inv{n} \sum_{k=0}^{n-1} 
\rho(X_k,\tilde X_{k+1})f(\tilde X_{k+1})+(1-\rho(X_k,\tilde
X_{k+1}))f(X_k).
\end{align*}
We shall  study in Section  \ref{sec:changingK} another choice  for the
weights also  considered by Frenkel  \cite{f:wrmc}.  Notice that  the WR
algorithm requires the  evaluation of $f$ for all  the proposals whereas
the Metropolis-Hastings  algorithm evaluates  $f$ only for  the accepted
proposals.   Other  algorithms using  all  the  proposals,  such as  the
Rao-Blackwell Metropolis-Hasting  algorithm, have been  studied, see for
example section 6.4.2 in  \cite{cr:mcsm} and references therein.  In the
Rao-Blackwell  Metropolis-Hasting  algorithm,  the weight  of  $f(\tilde
X_{k+1})$ depends on all the  proposals $\tilde X_1, \ldots, \tilde
X_n$. It is thus necessary to keep in memory the values of all
proposals in order to compute the estimation of $\langle \pi, f
\rangle$. 

One   easily   checks   that   $I^{WR}_n(f)-I_n(f)=J_n(f)$   where   for
any real function $\psi$ defined on $E$, 
\begin{align*}
   J_n(\psi)
&=\inv{n} \sum_{k=0}^{n-1} 
\bigg( \E\left[\psi(X_{k+1})|X_k, \tilde X_{k+1}\right]-\psi(X_{k+1})\bigg) \\
&=\inv{n} \sum_{k=0}^{n-1}\bigg( \rho(X_k,\tilde
     X_{k+1})\psi(\tilde
   X_{k+1})+(1-\rho(X_k,\tilde
  X_{k+1}))\psi(X_{k})-\psi(X_{k+1})\bigg) .
\end{align*}
Notice that $J_n(\psi)=0$ when $\psi$ is constant. We can consider a
more general estimator of $\langle \pi, f \rangle$ given by 
\[
I_n(f,\psi)=I_n(f)+J_n(\psi).
\]
Notice that $I^{WR}_n(f)=I_n(f,f) $ and $I_n(f)=I_n(f,0)$. It is easy to
check that the  bias of the estimator $I_n(f,\psi) $  does not depend on
$\psi$: $\E[I_n(f,\psi)]=\E[I_n(f)]$.  Theorem \ref{theo:princ2} implies
the following result on the estimator $I_n(f,\psi)$.
\begin{prop}
  \label{prop:CV}
   For any  real functions $\psi$ and $f$  defined on $E$,
  the  estimator   $I_n(f,\psi)  $  of  $\langle  \pi,   f  \rangle$  is
  consistent:   a.s.  $\displaystyle   \lim   _{n\rightarrow  \infty   }
  I_n(f,\psi)=\langle \pi, f \rangle$.

 \end{prop}

 {F}rom this result, $J_n(\psi)$ can be seen as a control variate and it
 is  natural to  look for  $\psi$ which  minimizes the  variance  or the
 asymptotic  variance of  $I_n(f,\psi)$.  Another class of control
 variates has been studied in \cite{atchperr} in the particular case of
 the Independent Metropolis-Hastings algorithm where $Q(x,.)$ does not
 depend on $x$. 

The last part of Theorem
 \ref{theo:princ2} implies the following result, where we used Lemma
 \ref{lem:added} to derive  the asymptotic variance expression. We shall
 write  $\E_\pi$ when $X_0$ is distributed under its invariant measure
 $\pi$ (in particular $\langle \pi, f \rangle=\E_\pi [f(X_0)]$). 
\begin{prop}\label{asvar} 
 For any real functions $\psi$ and $f$ defined on $E$, the estimator
$I_n(f,\psi) $ of $\langle \pi, f \rangle$  is asymptotically normal: 
$$\sqrt{n} \left(I_n(f,\psi)- \langle \pi,f \rangle\right) 
\; \xrightarrow[n\rightarrow \infty ]{(d)} \;
\cn(0,\sigma(f,\psi)^2),$$
with asymptotic variance $\sigma(f, \psi)^2$ given by
\begin{multline*}
\sigma(f,\psi)=\sigma(f)^2 -
 \E_\pi\left[\Big(1-\rho(X_0,
  X_1)\Big)\Big(F(X_1)-F(X_0)\Big)^2\right]\\
+ \E_\pi\left[\Big(1-\rho(X_0,
  X_1)\Big)\Big(\psi(X_1) -F(X_1)-\psi(X_0)+F(X_0)\Big)^2\right],  
\end{multline*}
where $F$ solves the Poisson equation \reff{Poisson}. In particular, for
fixed  $f$, the  asymptotic variance  $\sigma(f,\psi)^2$ is  minimal for
$\psi=F$ and  this  choice  achieves  variance reduction  :
$\sigma(f,F)^2\leq \sigma(f)^2$.
\end{prop}
Although  optimal in  terms of  the asymptotic  variance,  the estimator
$I_n(f,F)$ is not for use in practice, since computing a solution of the
Poisson  equation  is more  complicated  than  computing $\langle  \pi,f
\rangle$.    Nevertheless,   the   Proposition  suggests   that   using
$I_n(f,\psi)$ where $\psi$  is an approximation of $F$  might lead to a
smaller  asymptotic variance  than in  the  standard Metropolis-Hastings
algorithm.  Some hint at
the computation  of  an approximation of $F$ by a Monte  Carlo approach is
for  instance given in  \cite{munos} p.418-419.   Because of  the series
expansion $F=\sum_{k\geq 0}P^k(f-\langle \pi,f\rangle)$, $f$ can be seen
as an approximation of $F$  of order $0$.  Hence the asymptotic variance
of $I^{WR}_n(f)=I_n(f,f)$ might be smaller  than the one of $I_n(f)$ in
some situations.  It is  common belief in  the physicist  community, see
\cite{cck} or \cite{f:wrmc}, that  the inequality is always true. Notice
that, as  remarked by  Frenkel in a  particular case  \cite{f:wrmc}, the
variance of  each term of the  sum in $I_n^{WR}(f)$ is  equal or smaller
than the variance of each term  of the sum in $I_n(f)$ by Jensen inequality. 
But one has also to compare  the covariance terms,
which is not so obvious. We investigate whether the asymptotic
variance of the WR algorithm is smaller than the one of the
standard Metropolis algorithm and reach the following conclusion which
contradicts the intuition. 

\begin{prop}\label{prop:MB}$ $

\begin{itemize}
   \item[i)] In the Metropolis case, that is when \eqref{casga} holds with $\gamma(u)=\min(1,u)$, then it
     may happen that $\sigma(f,f)^2>\sigma(f)^2$.
\item[ii)] When \eqref{casga}
  holds with $\displaystyle \gamma(u)=\frac{\alpha u}{1+u}$, for some
  $\alpha\in (0,2)$, then we have
  $\sigma(f,f)^2\leq 
  \sigma(f)^2$.  Furthermore, for $f$ non constant, the function
  $b\mapsto \sigma(f, bf)^2$ is minimal at 
\begin{equation}
   \label{eq:bstar}
b_\star=\frac{\langle \pi, f^2 \rangle -\langle \pi, f \rangle^2}{\langle
  \pi, f^2 -fPf \rangle} 
\end{equation}
and $b_\star\geq 1/\alpha$. When $\alpha=1$, if, moreover, $\sigma(f,f)^2>0$, then $b_\star>1$.
\end{itemize}
\label{metbark}\end{prop}
\begin{rem}
   \label{rem:added}
Assume that $f$ is not constant. The optimal parameter $b_\star$ defined by
\reff{eq:bstar} can be estimated by 
\[
\hat b_n=\frac{I_n(f^2) -I_n(f)^2}{I_n(f^2) - \inv{n}\sum_{k=1}^n
  f(X_{k-1})f(X_k)}\cdot
\]
Notice that a.s. $\lim_{n\rightarrow \infty } \hat b_n=b_\star$ thanks to
the ergodic theorem. Using Slutsky theorem, one can deduce from
Proposition \ref{asvar}  that $I_n(f)+ \hat b_n J_n(f)=I_n(f, \hat b_n
f) $ is an asymptotically normal estimator of $\langle \pi, f \rangle$
with asymptotic variance $\sigma(f, b_\star f)^2$. Thus, in the framework ii)
of Proposition \ref{prop:MB}, using the control variate $\hat b_n
J_n(f)$ improves strictly the WR estimator as soon as either $\alpha<1$ or
$\alpha=1$ (Boltzmann algorithm) and $\sigma(f,f)^2$ is positive. Notice that when
 its asymptotic
variance $\sigma(f,f)^2$ is zero, then the WR estimator $I_n^{WR}(f)=I_n(f,f)$ is equal
to $\langle \pi,f \rangle$.
\end{rem}

To prove assertion i), we give an explicit counter-example such that
$\sigma(f,f)^2>\sigma(f)^2$ in the Metropolis case (see Section 
\ref{sec:ce} and equation  \reff{eq:sff>sf}). The assertion ii) is
also proved in Section \ref{sp} (see Proposition \ref{prop:geneB}). Let us
make some comments on its hypothesis which holds with $\alpha=1$ for
Boltzmann acceptation rule.  
\begin{itemize}
\item By \eqref{prerevers} and since $\rho(x,y)$ is an acceptance probability, the constant $\alpha $  has to be
smaller than 
$\displaystyle 1+\min_{x\neq   y,
  Q(x,y)>0}\frac{\pi(y)Q(y,x)}{\pi(x)Q(x,y)}$. 
\item If there exists a constant $c>0$ s.t. 
for
all distinct $x,y\in E$ s.t. $Q(x,y)>0$, the quantity 
$\displaystyle \frac{\pi(x)Q(x,y)}{\pi(y) Q(y,x)}$ is  equal to $c$ or
$1/c$ and \eqref{casga} holds with $\gamma$ such that
$\gamma(1/c)=\gamma(c)/c$ then the hypothesis holds with $\alpha=\gamma(c)+\gamma(1/c)$.  For example assume that the transition
matrix $Q$ is symmetric and that $\pi$ is written as a Gibbs
distribution: for all $x\in E$, $ \pi(x) = 
\expp{- H(x)}/\sum_{y\in
  E} \expp{-H(y)}$ for some energy function $H$. If  the energy increases or
decreases by the 
same amount $\varepsilon$ for all the authorized transitions, then
$\displaystyle \frac{\pi(x)Q(x,y)}{\pi(y) Q(y,x)}$ is equal to $c$ or
$1/c$ with $c=e^{\varepsilon}$.
\end{itemize}

According  to \cite{p:omcsumc},  since for  all  $u>0$, $\displaystyle
  \frac{u}{1+u}<\min(1,u)$,  in  the  absence  of waste  recycling,  the
  asymptotic variance  $\sigma(f)^2$ is  smaller in the  Metropolis case
  than  in the  Boltzmann case  for given  $\pi$, $Q$  and $f$.  So waste
  recycling always achieves variance reduction only for the worst choice
  of $\gamma$. Notice however that the Boltzmann algorithm is used in
  the multi-proposal framework where we generalize our results. 

\begin{rem}
\label{rem:Reduc}
   When the computation of $Pg$ is feasible for any function
   $g:E\rightarrow\R$ (typically when, for every $x\in E$, the cardinal of $\{y\in
   E:Q(x,y)>0\}$ is small), then it is possible to use $I_n(\psi-P\psi)$
   as a control variate and approximate $\langle \pi,f\rangle$ by
   $I_n(f-(\psi-P\psi))$. Since $\pi$ is invariant with respect to $P$,
   $\langle \pi,\psi-P\psi\rangle=0$ and a.s. $I_n(f-(\psi-P\psi))$
   converges to $\langle \pi,f\rangle$ as $n$ tends to
   infinity. Moreover, the asymptotic variance of the estimator is
   $\sigma(f-\psi+P\psi)^2$. Last, remarking that 
\begin{equation}
   I_n(\psi-P\psi)=\frac{1}{n}\sum_{k=1}^n\left(\psi(X_k)-P\psi(X_{k-1})\right)+\frac{1}{n}\left(P\psi(X_0)-P\psi(X_n)\right)\label{inppp}
\end{equation}
one obtains that the bias difference 
$\E[I_n(f-\psi+P\psi)]-\E[I_n(f)]=\frac{1}{n}\E\left[P\psi(X_0)-P\psi(X_n)\right]$
is smaller than $2\max_{x\in E}|\psi(x)|/n$.

For the choice $\psi=F$, this control variate is perfect, since
according to \eqref{Poisson}, for each $n\in\N^*$, $I_n(f-(F-PF))$ is
constant and equal to $\langle \pi,f\rangle$.

For the choice $\psi=f$, the asymptotic
variance of the estimator $I_n(Pf)$ is also smaller than the one of $I_n(f)$. Indeed setting
$f_0=f-\langle \pi,f\rangle$, we have
\begin{align*}
   \sigma(f)^2-\sigma(Pf)^2
&=\langle \pi,F^2 + (P^2 F)^2 - 2 (PF)^2\rangle\\
&=\langle \pi,(f_0+PF)^2-2(PF)^2+(Pf_0-PF)^2\rangle\\
&=\langle \pi,f_0^2+2f_0P(F-PF)+(Pf_0)^2\rangle=\langle \pi,(f_0+Pf_0)^2\rangle
\end{align*}
where we used that $PF$ solves the Poisson equation \eqref{Poisson} with
$f$  replaced  by $Pf$  and  \reff{eq:def-s2}  for  the first  equality,
\eqref{Poisson} for the second and last equalities and the reversibility
of $\pi$ w.r.t. $P$ for the last one.

Notice the control variate $J_n(\psi)$ is similar to $I_n(\psi-P\psi)$
except that 
the conditional expectation $P\psi(X_{k-1})$ of $\psi(X_k)$ given
$X_{k-1}$ in the first term of the r.h.s. of \eqref{inppp} is replaced by the
conditional expectation of $\psi(X_k)$ given
$(X_{k-1},\tilde{X}_k)$ which can always be easily computed. From this
perspective, the minimality of the asymptotic variance of 
$I_n(f,\psi)$ for $\psi=F$ is not a surprise.

The comparison between $\sigma(f,\psi)^2$ and $\sigma(f-\psi+P\psi)^2$
can be deduced from Section \ref{sec:compvarppp} which is stated in the more
general multi-proposal framework introduced in the next
paragraph. Notice that 
the sign of $\sigma(f,\psi)^2-\sigma(f-\psi+P\psi)^2$ depends on $\psi$.
\end{rem}
\subsection{Multi-proposal WR algorithm}
In  the  classical  Metropolis  Hasting  algorithm, there  is  only  one
proposal $\tilde{X}_{n+1}$  at step $n+1$. Around  1990, some extensions
where  only one  state among  multiple proposals  is accepted  have been
proposed  in order to  speed up  the exploration  of $E$ (see
\cite{ad:hrud} 
for a unifying presentation of MCMC algorithms including the
multi-proposal Metropolis  Hasting  algorithm).   According to
Frenkel \cite{f1:wrmc},  the waste recycling can  be particularly useful
for these algorithms where many states are rejected.

To formalize
these     algorithms,    we     introduce    a     proposition    kernel
$\cQ:E\times{\mathcal P}(E)\rightarrow  [0,1]$, where ${\mathcal P}(E)$ denotes the set of parts
of $E$, which describes how to randomly choose the set of proposals: 
\begin{equation}
   \label{eq:propQ}
\forall x\in
E,\;\cQ(x,A)=0\mbox{ if }x\notin A\quad \mbox{ and } \quad \sum_{A\in
  {\mathcal P}(E)}\cQ(x,A)=1.
\end{equation}
The second condition says that $ \cQ(x,\cdot)$ is a probability on
$\cp(E)$. The first one ensures that the starting point is among the
proposals. This last convention will allow us to transform the
rejection/acceptation procedure  into a selection procedure among the
proposals.  

The  selection procedure is  described by  a probability  $\kappa$.  For
$(x,A)\in E\times{\mathcal  P}(E)$, let $\kappa(x,A,\tilde{x})\in [0,1]$
denote the  probability of choosing  $\tilde{x}\in A$ as the  next state
when  the   proposal  set   $A$  has  been   chosen.   We   assume  that
$\sum_{\tilde{x}\in  A}\kappa(x,A,\tilde{x})=1$ (that  is  $\kappa(x, A,
\cdot)$ is a probability measure) and that the following condition holds :
\begin{equation}
   \forall A\in{\mathcal P}(E),\;\forall x,\tilde{x}\in
   A,\;\pi(x)\cQ(x,A)\kappa(x,A,\tilde{x})
   =\pi(\tilde{x})\cQ(\tilde{x},A)\kappa(\tilde{x},A,x).  
\label{reversmult}
\end{equation}
This condition is the  analogue of \reff{prerevers} for a multi-proposal
setting. For examples of non-trivial selection probability $\kappa$, see
after Proposition \ref{prointmul}.

The Markov  chain $X=(X_n,n\geq  0)$ is now  defined inductively  in the
following way.  Let $X_0$ be a random variable taking values in $E$ with
probability distribution  $\nu_0$. At step  $n$, $X_0, \ldots,  X_n$ are
given.  The  proposal  set  at  step $n+1$,  $A_{n+1}$,  is  distributed
according  to  $\cQ(X_n,\cdot)$. Then  $X_{n+1}$  is chosen  distributed
according to $\kappa(X_n,A_{n+1},.)$. It is  easy to check that $X$ is a
Markov chain with transition matrix
\begin{equation}
   P(x,y)=\sum_{A\in{\mathcal P}(E):x,y\in A} 
\cQ(x,A)\kappa(x,A,y).\label{defPmult}
\end{equation}
Condition
\eqref{reversmult} ensures that $X$ is reversible w.r.t. the probability
measure $\pi$ : $\pi(x)P(x,y)=\pi(y)P(y,x)$. 
\begin{rem}
   \label{rem:multi-1}
The multi-proposal Metropolis-Hastings algorithm 
generalizes the Metropolis-Hastings algorithm  which can be
recovered for the particular 
choice $\cQ(x,\{x,y\})=Q(x,y)$ and for $y\neq x$,
$\kappa(x,\{x,y\},y)=1-\kappa(x,\{x,y\},x)=\rho(x,y)$. 
\end{rem}

We keep the definition \reff{eq:def-In} of $I_n(f)$ but adapt the ones
of $J_n(\psi)$ and 
$I_n(f,\psi)$ as follows : 
\begin{align}
\nonumber
  \cJ_n(\psi)
&=\inv{n} \sum_{k=0}^{n-1}\bigg(\E\left[\psi(X_{k+1})|X_k,A_{k+1}\right]
-\psi(X_{k+1})\bigg) \\
&=\inv{n} \sum_{k=0}^{n-1}\bigg(\sum_{\tilde{x}\in
  A_{k+1}}\kappa(X_k,A_{k+1},\tilde{x})\psi(\tilde{x})-\psi(X_{k+1})\bigg)
\label{cjpsi}  
\end{align}
and  $\cI_n(f,\psi)=I_n(f)+\cJ_n(\psi)$.  The Waste  Recycling estimator
of $\langle \pi,f\rangle$ studied  by Frenkel in \cite{f1:wrmc} is given
by  $\cI^{WR}_n(f)=\cI_n(f,f)$.   Notice  that the  bias   of  the  estimator
$\cI_n(f,\psi)     $    does    not     depend    on     $\psi$    (i.e.
$\E[\cI_n(f,\psi)]=\E[I_n(f)]$).
  It   turns  out   that   Propositions
\ref{prop:CV}  and  \ref{asvar}   remain  true  in  this  multi-proposal
framework (see Theorem \ref{theo:princ2}) as soon as $P$ is irreducible.
Notice that the irreducibility of $P$ holds  if and only if for all $ x'\neq
y\in E$, there exist $m\geq 1$, distinct $x_0=y,x_1,x_2,\hdots,x_m=x'\in
E$   and  $A_1,A_k\hdots,A_m\in{\mathcal  P}(E)$   such  that   for  all
$k\in\{1,\hdots,m\}$, $x_{k-1},x_k\in A_k$ and
\begin{equation}
\label{irregen}
   \prod_{k=1}^m \cQ(x_{k-1},A_k) \kappa(x_{k-1}, A_k,x_k)>0. 
\end{equation}
\begin{prop}\label{prointmul}
Assume that $P$ is irreducible. For any real functions $\psi$ and $f$ defined on $E$, we have:
\begin{itemize}
   \item The estimator
$\cI_n(f,\psi) $ of $\langle \pi, f \rangle$  is  consistent:
a.s. $\displaystyle  \lim _{n\rightarrow \infty } \cI_n(f,\psi)=\langle
\pi, f \rangle$.   
    \item The estimator
$\cI_n(f,\psi) $ of $\langle \pi, f \rangle$  is asymptotically normal: 
$$\sqrt{n} \left(\cI_n(f,\psi)- \langle \pi,f \rangle\right) 
\; \xrightarrow[n\rightarrow \infty ]{(d)} \;
\cn(0,\sigma(f,\psi)^2)$$
where the asymptotic variance (still denoted by) $\sigma(f, \psi)^2$ is
given by 
\[
\sigma(f,\psi)^2=\sigma(f)^2 + \sum_{x\in E, A\in \cp(E)} \pi(x)
\cQ(x,A) \left[\Var_{\kappa_{x,A}}(\psi-F)-\Var_{\kappa_{x,A}}(F)
\right],
\]
with  $\displaystyle \Var_{\kappa_{x,A}}(g)=\sum_{y\in A} \kappa(x,A,y)
g(y)^2 - \left(\sum_{y\in A} \kappa(x,A,y)
g(y)\right)^2$. 
   \item Moreover, for fixed $f$, the asymptotic variance
$\sigma(f,\psi)^2$ is minimal for $\psi=F$ where $F$ solves the Poisson
equation \reff{Poisson}. In particular, this choice achieves
variance reduction: $\sigma(f,F)^2\leq
\sigma(f)^2$.
\end{itemize}
\end{prop}

We now give two examples of non-trivial
selection probability $\kappa$ which satisfies condition
\reff{reversmult}. The first one, $\kappa^M$, 
defined by 
\begin{equation}
\label{defkapm}\kappa^{M}(x,A,\tilde{x})=\begin{cases}\displaystyle
  \frac{\pi(\tilde{x})\cQ(\tilde{x},A)}{\max
    \left(\pi(\tilde{x})\cQ(\tilde{x},A), \pi(x)\cQ(x,A)\right)+\sum_{z\in
    A\setminus \{x,\tilde{x}\}} \pi(z)\cQ(z,A)}&\quad
\text{if}\quad\tilde{x}\neq x,\\ 
1-\sum_{z\in A\setminus \{x\}}\kappa^M(x,A,z)&\quad
\text{if}\quad\tilde{x}=x, 
\end{cases}
\end{equation}
generalizes the
Metropolis selection given by \eqref{casga} with $\gamma(u)=\min(1,u)$. 
(Notice that for $x\neq \tilde x$ one has $\displaystyle \kappa^{M}(x,A,\tilde{x})
\leq  \frac{\pi(\tilde{x})\cQ(\tilde{x},A)}{\sum_{z\in
    A\setminus \{x\}} \pi(z)\cQ(z,A)}$, which implies that $1-\sum_{z\in
  A\setminus \{x\}}\kappa^M(x,A,z)$ is indeed non-negative.) 
The second one,  $\kappa^B$, which does not depend on the initial point $x$,
and is defined by 
\begin{equation}
\label{defkapb}
   \kappa^{B}(x,A,\tilde{x})=\kappa^{B}(A,\tilde{x})
   =\frac{\pi(\tilde{x})\cQ(\tilde{x},A)}{\sum_{z\in  
    A} \pi(z)\cQ(z,A)},
\end{equation}
generalizes the  Boltzmann (or Barker) selection  given by \eqref{casga}
with  $\displaystyle \gamma(u)=\frac{u}{1+u}  $.  Notice  that  for both
choices, the  irreducibility condition \eqref{irregen}  can be expressed
only in terms of $\cQ$ :
\begin{equation*}
   \prod_{k=1}^m\cQ(x_{k-1},A_k)\cQ(x_{k},A_k)>0.
\end{equation*}

For  the  selection  probability  \reff{defkapb}, we  prove  in  section
\ref{sec:Bolt}  (see   Proposition  \ref{propbarker})  that   the  Waste
Recycling improves the Metropolis-Hasting algorithm :
\begin{prop}
\label{prop:mult-MB}
When $\kappa=\kappa^B$ is given by \eqref{defkapb} (Boltzmann
  or Barker case), then we have $\sigma(f,f)^2\leq
  \sigma(f)^2$. Furthermore, for $f$ non constant, the function
  $b\mapsto \sigma(f, bf)^2$ is minimal at $b_\star$ defined by
  \reff{eq:bstar} and $b_\star>1$ when $\sigma(f,f)^2>0$. 
\end{prop}
Since   for   $\tilde{x}\neq   x\in   A$,   $\kappa^M(x,A,\tilde{x})\geq
\kappa^B(A,\tilde{x})$,  according to  \cite{p:omcsumc},  the asymptotic
variance $\sigma(f)^2$  remains smaller in  the Metropolis case  than in
the  Boltzmann one.  Nethertheless,  it is  likely  that the  difference
decreases when  the cardinality of  the proposal sets  increases. Notice
that the optimal value $b_\star$ can be estimated by $\hat b_n$ which is computed
using  the proposals: see Remark  \ref{rem:added}. The  control variate
$\hat b_n J_n(f)$ improves therefore the WR algorithm.

\section{Main result for general multi-proposal WR} 
\label{sec:Egene}

Let  $(E,  \cf_E)$  be  a  measurable  space s.t. $\{x\}\in \cf_E$ for
all $x\in E$, and  $\pi$  be  a  
probability measure    on    $E$. Notice that $E$ is not 
assumed to be 
finite. Let $\displaystyle
\cp=\{A\subset E; \Card(A)<\infty \}$ be the set of finite subsets of
$E$. Let $\displaystyle \bar E= \cup_{n\geq 1} E^n$ and $\cf_{\bar E}$
the smallest 
$\sigma$-field on $\bar E$ which contains $A_1\times \cdots
\times A_n$ for all $A_i \in \cf_E$ and $n\geq 1$. 
We consider the function $\Gamma$ defined on $\bar E $ taking
value on $\cp$ such that $\Gamma ((x_1, \ldots, x_n))$ is the set
$\{x_1, \ldots, x_n\}$ of distinct
elements in $(x_1, \ldots, x_n)$. We define $\cf_\cp$, a
$\sigma$-field on $\cp$, as the image of $\cf_{\bar E} $ by the
application $\Gamma$.  
We consider a measurable proposition probability kernel 
$\cQ:E\times \cf_\cp\rightarrow   [0,1]$  s.t.  
\begin{equation}
   \label{eq:propQ2}
\int_\cp \cQ(x, dA)=1\quad\text{and}\quad \int_\cp \cQ(x, dA)\; \ind_{\{
  x\not \in A\}}=0
\end{equation}
(this  is the analogue  of \reff{eq:propQ})  and a  measurable selection
probability kernel $\kappa: E\times  \cp \times \cf_E \rightarrow [0,1]$
s.t.   for $x\in  A$ we  have $\kappa(x,A,A)=1$.  Let $\delta_y$  be the
Dirac mass  at point  $y$. In  particular, since $A$  is finite,  with a
slight  abuse of notation,  we shall  also write  $\kappa(x,A,dy)
=\sum_{z\in A} \kappa(x,A,z) \delta_z(dy)$ and so $\sum_{y\in A}
\kappa(x,A,y) =1$.

We assume that the analogue of \eqref{reversmult} holds, that is 
\begin{equation}
   \label{eq:reversmult2}
\pi(dx)\cQ(x,dA)\kappa(x,A,dy)
   =\pi(dy)\cQ(y,dA)\kappa(y,A,dx).  
\end{equation}

\begin{exa}
\label{exa:BMgen}
We give  the analogue of  the Metropolis and Boltzmann   selection kernel
defined  in \reff{defkapm} and  \reff{defkapb} when  $E$ is  finite.  We
consider   $N(dx,  dA)=\pi(dx) \cQ(x,dA)$  and  a  measure
$N_0(dA)$  on  $\cf_\cp$  such  that  $\int_{x\in  E}  N(dx,  dA)  $  is
absolutely continuous w.r.t. $N_0(dA)$.  Since $x\in A$ and $A$ is finite
$N(dx,  dA)$-a.s., the  decomposition  of $N$  w.r.t.  $N_0$ gives  that
$N(dx,  dA)=N_0(dA)   r_A(dx)$,  where  $r_A(dx)=\sum_{y\in   A}  r_A(y)
\delta_y(dx)$ if $A$ is finite and $r_A(dx)=0$ otherwise, and $(x,A)
\mapsto r_A(x)$ is jointly measurable.

The Metropolis selection kernel is given by: for $x,y\in A$, $r_A\neq
0$, 
\begin{equation}
   \label{eq:Mgen}
\kappa^M(x,A,y)=\frac{r_A(y)}{\sum_{z\in
    A\setminus \{x,y\}} r_A(z)+\max(r_A(x), r_A(y))},
\end{equation}
if $x\neq y$  and $\kappa^M(x,A,x)=1-\sum_{y\in
  A\setminus\{x\}}\kappa^M(x,  A,y)$.

The Boltzmann selection kernel is given by: for $x,y\in A$, $r_A\neq 0$, 
\begin{equation}
   \label{eq:Bgen}
\kappa^B(x,A,y)=\kappa^B(A,y)=\frac{r_A(y)}{\sum_{z\in A} r_A(z)}.
\end{equation}
We  choose  those two  selection  kernels to  be  equal  to the  uniform
distribution  on $A$  when $r_A=0$.   For those  two  selection kernels,
equation \reff{eq:reversmult2} is satisfied.
\end{exa}

\begin{exa}
\label{exa:Exple-cont}
Let us give a natural example.   Let $\nu$ be a reference measure on $E$
with no  atoms, $\pi$ a probability  measure on $E$  with density w.r.t.
$\nu$ which  we still  denote by $\pi$,  a selection procedure  given by
$\cQ(x,\ca)=\P_x(\{x,Y_1,  \ldots,  Y_n\}\subset   \ca)$  for  $\ca  \in
\cf_\cp$,  where $Y_1,  \ldots, Y_n$  are $E$-valued  independent random
variables with density w.r.t. $\nu$  given by $q(x, \cdot)$ under $\P_x$
and $n\geq  1$ is fixed.   We use notations of  Example \ref{exa:BMgen}.
In  this setting, we  choose $N_0(dA)=\prod_{x\in  A}  \nu(dx)$ and  the
function  $r_A$  is  given  by:  for $x\in  A$,  $\displaystyle  r_A(x)=
\pi(x)\prod_{z\in A\setminus \{x\} }q(x,z)$.
\end{exa}

The Markov  chain $X=(X_n,n\geq  0)$ is  defined inductively  in the
following way.  Let $X_0$ be a random variable taking values in $E$ with
probability distribution  $\nu_0$. At step  $n$, $X_0, \ldots,  X_n$ are
given.  The  proposal  set  at  step $n+1$,  $A_{n+1}$,  is  distributed
according  to  $\cQ(X_n,\cdot)$. Then  $X_{n+1}$  is chosen  distributed
according to $\kappa(X_n,A_{n+1},.)$. This is a particular case of the
hit and run algorithm \cite{ad:hrud}, where the proposal sets are always
finite. It is  easy to check that $X$ is a
Markov chain with transition kernel
\begin{equation}
   \label{eq:defP}
   P(x,dy)=\int_\cp \cQ(x,dA)\kappa(x,A,dy).
\end{equation}
For $f$ a real valued measurable function defined on $E$, we shall write
$Pf(x)$ for $\int_E P(x,dy) f(y)$ when this integral is well defined. 

Condition \eqref{eq:reversmult2}  ensures that $X$  is reversible w.r.t.
$\pi$  : $\pi(dx)P(x,dy)=\pi(dy)P(y,dx)$.  We  also assume  that $X$  is
Harris recurrent  (see \cite{mt:mcss} section 9). This  is equivalent to
assume  that for  all $B\in  \cf_E$  s.t. $\pi(B)>0$  we have  $\P(\Card
\{n\geq 0; X_n\in B\} =\infty|X_0=x )=1$  for all $x\in E$. 

\begin{exa}
It is easy to
check in  Example  \ref{exa:Exple-cont} that $X$ is
Harris recurrent if the random walk with transition kernel $q$ is itself
Harris recurrent and $$\forall x\in E,\;\cQ(x,dA)\mbox{ a.e. },\;\forall y\in
A,\;\kappa(x,A,y)>0.$$
\end{exa}

For $f$  a real valued  measurable function defined  on $E$ and  $\nu$ a
measure on $E$, we shall  write $\langle \nu,f\rangle$ for $\int \nu(dy)
f(y)$ when this integral is well defined.

Let  $f$  be a  real-valued  measurable  function  defined on  $E$  s.t.
$\langle \pi,  |f| \rangle <\infty $.  Theorem  17.3.2 in \cite{mt:mcss}
asserts that  a.s.  $\lim_{n\rightarrow  \infty } I_n(f)=\langle  \pi, f
\rangle$, with $I_n(f) $ defined by \reff{eq:def-In}.

We consider  the functional $\cJ_n$ defined by 
\begin{align}
\nonumber
   \cJ_n(\beta)
&=\inv{n} \sum_{k=0}^{n-1}\bigg( \E\left[\beta(X_k, A_{k+1}, X_{k+1})|X_k,
  A_{k+1}\right]   -\beta(X_k,A_{k+1},X_{k+1})\bigg)\\
&=\inv{n} \sum_{k=0}^{n-1}\bigg(\;\;\sum_{\tilde{x}\in
   A_{k+1}} \kappa(X_k,A_{k+1},\tilde{x})
   \beta(X_k,A_{k+1},\tilde{x})   -\beta(X_k,A_{k+1},X_{k+1})\bigg),
\label{cjbet} 
\end{align} 
for $\beta$ any real-valued  measurable function defined on $E\times \cp
\times E$.   We set $\cI_n(f,\beta)=I_n(f)+\cJ_n(\beta)$.   To prove the
convergence and the asymptotic normality of the estimator $\cI_n(f, \beta) $
of  $\langle \pi,  f\rangle$, we  shall  use a  martingale approach.  In
particular, we shall  assume there exists $F$ a  solution to the Poisson
equation $F-PF=f -\langle \pi, f \rangle$ s.t. $\langle \pi, F^2 \rangle
<\infty   $  (see   theorem  17.4.2   and  condition   (V.3)   p.341  in
\cite{mt:mcss} to ensure the existence of such a solution).

We  introduce  the  following  convenient notation.  For  a  probability
measure $\nu$  on $E$ and real  valued functions $h$ and  $g$ defined on
$E$, we write, when well defined,
\[
\Cov_\nu(h,g)=\langle \nu, gh \rangle -\langle \nu,g \rangle\langle
\nu,h \rangle \quad\text{and}\quad \Var_\nu(h)=\langle \nu, h^2 \rangle
-\langle \nu,h \rangle^2
\]
respectively the covariance of $g$ and $h$ and the variance of $h$
w.r.t. $\nu$.  We also write $\kappa_{x,A}(dy)$ for the probability measure
$\kappa(x,A,dy )$ and the $\beta_{x,A}(\cdot)$ for the function
$\beta(x,A,\cdot)$. 

\begin{theo}
   \label{theo:princ2}
We assume $X$ is Harris recurrent, $\langle \pi, f^2 \rangle<\infty $,
there exists a solution $F$ to  the Poisson equation $F-PF=f
-\langle  \pi, f \rangle$ such that $\langle
\pi, F^2 \rangle<\infty $,  and $\beta$ is square integrable:  $\int \pi(dx) \cQ(x,dA)
\kappa(x,A,dy) \beta(x,A,y)^2<\infty $. Under those  assumptions, we have:
\begin{itemize}
   \item[$(i)$] The estimator
$\cI_n(f,\beta) $ of $\langle \pi, f \rangle$  is  consistent:
a.s. $\displaystyle  \lim _{n\rightarrow \infty } \cI_n(f,\beta)=\langle
\pi, f \rangle$.   
    \item[$(ii)$] The estimator
$\cI_n(f,\beta) $ of $\langle \pi, f \rangle$  is asymptotically normal:
\[
\sqrt{n} \left(I_n(f,\beta)- \langle \pi,f \rangle\right) 
\;        \xrightarrow[n\rightarrow        \infty       ]{(d)}        \;
\cn(0,\sigma(f,\beta)^2),
\]
and the  asymptotic variance is given by 
\begin{equation}
   \label{varFdel}
\sigma(f,\beta)^2=\sigma(f)^2+\int
 \pi(dx)\cQ(x,dA)\left[\Var_{\kappa_{x,A}} (\beta_{x,A}-F) 
  - \Var_{\kappa_{x,A}} (F)\right],
\end{equation}
with $\sigma(f)^2=\langle \pi, F^2 -(PF)^2 \rangle$. 
   \item[$(iii)$]   The asymptotic variance $\sigma(f,\beta)^2$ is minimal for
    $\beta_{x,A}=F$  and 
\begin{equation}
   \label{eq:sfF2}
\sigma(f,F)^2=\int \pi(dx)\left(\int\cQ(x,dA)\langle \kappa_{x,A}, F
  \rangle^2 - \left(\int\cQ(x,dA)\langle \kappa_{x,A}, F
    \rangle\right)^2\right)\leq 
 \sigma(f)^2.
\end{equation}
\end{itemize}
\end{theo}

\begin{proof}
We shall prove the Theorem when $X_0$ is  distributed according to
$\pi$. The general case follows from proposition 17.1.6 in \cite{mt:mcss}, 
since $X$ is Harris
recurrent.

  We set, for $n\geq 1$,
\[
\Delta M_{n}=F(X_{n})-PF(X_{n-1}) + \eta (X_{n-1}, A_{n},X_{n}), 
\]
where  
\[
\eta (x,A,y)=  \sum_{\tilde{x}\in A}\left(\kappa(x,A,\tilde{x})-
  \ind_{\{y=\tilde{x}\}}\right)\beta(x,A,\tilde{x}).
\]
Notice that  $\Delta M_n$ is  square integrable and  that $\displaystyle
\E[\Delta  M_{n+1}|\cg_n]=0$,  where  $  \cg_n$  is  the  $\sigma$-field
generated by $X_0$ and $(A_i,  X_i)$ for $1\leq i\leq n$.  In particular
$M=(M_n, n\geq  0)$ with $M_n=\sum_{k=1}^n  \Delta M_k$ is  a martingale
w.r.t.  to the filtration $(\cg_n, n\geq 0)$.  Using that $F$ solves the
Poisson equation, we also have 
\begin{equation}
   \label{eq:In-M}
\cI_n(f,\beta)=\inv{n}M_n -\inv{n}PF(X_{n})+\inv{n}PF(X_{0})+ \langle
\pi, f \rangle.
\end{equation}
As $\langle \pi,  F^2 \rangle <\infty $ implies  that $\langle \pi, |PF|
\rangle <\infty $, we deduce  from theorem 17.3.3 in \cite{mt:mcss} that
a.s.  $\lim_{n\rightarrow \infty } \inv{n}PF(X_{n})=0$.  In particular
part  $(i)$  of  the  Theorem  will be  proved  as soon  as  we check  that
a.s. 
$\lim_{n\rightarrow \infty } \inv{n}  M_n=0$.  

We easily   compute
the bracket  of $M_n$:
\[
\langle M \rangle_n=\sum_{k=1}^n \E[\Delta  M_k^2|\cg_{k-1}]
=\sum_{k=1}^n h(X_{k-1}),
\]
with 
\[
h(x)=P(F^2)(x)- (PF(x))^2  + \int  \cQ(x,
dA)\left[-2\Cov_{\kappa(x,A,\cdot)} (\beta_{x,A}, 
     F) 
  + \Var_{\kappa(x,A,\cdot)} (\beta_{x,A})\right].
\]
Elementary computation yields 
\[
-2\Cov_{\kappa(x,A,\cdot)} (\beta_{x,A}, 
     F)  
  + \Var_{\kappa(x,A,\cdot)} (\beta_{x,A})= \Var_{\kappa(x,A,\cdot)}
  (\beta_{x,A}-F)  
  - \Var_{\kappa(x,A,\cdot)} (F).
\]
Since $\langle \pi,  F^2 \rangle <\infty $ and  $\int \pi(dx) \cQ(x, dA)
\kappa(x,  A, dy)  \beta(x,A,y)^2<\infty $,  we have  that $h$  is $\pi$
integrable. We  set $\sigma(f,\beta)^2=\langle \pi, h  \rangle$, that is
$\sigma(f,\beta)^2$    is   given    by   \reff{varFdel},    thanks   to
\reff{eq:def-s2} and the fact that  $\pi$ is invariant for $P$.  Theorem
17.3.2 in \cite{mt:mcss} asserts that a.s.  $\lim_{n\rightarrow \infty }
\inv{n}  \langle M  \rangle_n  =\langle \pi,  h  \rangle$. Then  theorem
1.3.15 in  \cite{d:rim} implies  that a.s. $\lim_{n\rightarrow  \infty }
\inv{n} M _n =0$. This ends the proof of part $(i)$.

The proof of part $(ii)$ relies on the central limit theorem for 
martingales, see
theorem 2.1.9 in \cite{d:rim}. We have already proved that a.s. 
 $\lim_{n\rightarrow  \infty } \inv{n}
\langle M \rangle_n =\sigma(f,\beta)^2$. Let us now check the
Lindeberg's condition. Notice that theorem  17.3.2 in
\cite{mt:mcss} implies that for any $a>0$, we have 
\[
\lim_{n\rightarrow\infty } 
\inv{n} \sum_{k=1}^n \E[\Delta M_k^2\ind_{\{|\Delta M_k^2>a\}}
|\cg_{k-1}] = \langle \pi, h_a \rangle, 
\]
where $h_a(x)=\E[\Delta M_1^2\ind_{\{|\Delta M_1^2>a\}}
|X_0=x]$. Notice that $0\leq h_a\leq h$ and that $(h_a, a>0)$ decreases
to $0$ as $a$ goes to infinity. We deduce that a.s. 
\[
\limsup_{n\rightarrow\infty } 
\inv{n} \sum_{k=1}^n \E[\Delta M_k^2\ind_{\{|\Delta M_k^2>\sqrt{n} \}}
|\cg_{k-1}] \leq \limsup_{a\rightarrow \infty } \langle \pi, h_a \rangle=0.
\]
This  gives the Lindeberg's condition. We deduce then that
$(\inv{\sqrt{n}}M_n , n\geq 1)$ converges in distribution to $\cn(0,
\sigma(f,\beta)^2)$. Then use \reff{eq:In-M} and that
a.s. $\displaystyle \lim_{n\rightarrow \infty } \inv{n}
(PF(X_{n+1}))^2=0$ 
(thanks to theorem 17.3.3 in \cite{mt:mcss}) 
to get  part $(ii)$. 

Proof of  part $(iii)$.  The asymptotic variance  $\sigma(f,\beta)^2$ is
minimal when $\Var_{\kappa_{x,A}} (\beta_{x,A}-F) =0$ that is at least
for $\beta_{x,A}=F$. Of course, $\sigma(f,F)^2\leq
 \sigma(f,0)^2=\sigma(f)^2$. 
Using \reff{eq:def-s2}, that $ \pi$ is invariant for $P$ and the
definition \reff{eq:defP} of $P$, we get 
\begin{align*}
   \sigma(f)^2
&=\langle \pi,PF^2\rangle-\langle \pi,(PF)^2\rangle\\
&= \int \pi(dx)\cQ(x,dA)
\langle \kappa_{x,A} , F^2  \rangle -\int
\pi(dx)\left(\int\cQ(x,dA)\langle \kappa_{x,A},F \rangle\right)^2.
\end{align*}
And the expression of $\sigma(f,F)^2$ follows from \eqref{varFdel}.

\end{proof}

\section{The Boltzmann case}
\label{sec:Bolt}
We  work in  the general  setting  of Section  \ref{sec:Egene} with  the
Boltzmann selection  kernel $\kappa$  given by \reff{eq:Bgen}  (or simply
\reff{defkapb}  when $E$  is finite).   The next  Proposition generalizes
Proposition \ref{prop:mult-MB}.   It ensures that  the asymptotic variance
of the waste recycling algorithm $\sigma(f,f)^2$ is smaller than the one
$\sigma(f)^2$  of the  standard Metropolis  Hastings algorithm and that
$b\mapsto \sigma(f,bf)^2$ is minimal at $b_\star$ given by \eqref{eq:bstar}.   In the
same time, we show that  this variance $\sigma(f)^2$ is at least divided
by two for the optimal choice $\beta(x,A,y)=F(y)$ in our control variate
approach.

For $f$ s.t. $\langle \pi, f^2 \rangle<\infty$, we set  $f_0=f-\langle
\pi,  f  \rangle$ and  
\begin{equation}
   \label{eq:Deltaf}
 \Delta(f)=\inv{2}\int     \pi(dx)    P(x,dy)
(f_0(x)+f_0(y))^2  = \langle \pi, f_0(f_0+Pf_0) \rangle. 
\end{equation}
Notice that the second equality in \reff{eq:Deltaf} is a consequence of the
invariance of $\pi$ w.r.t. $P$. 
\begin{prop}
\label{propbarker}
We assume that $X$ is Harris recurrent, $\langle \pi, f^2 \rangle<\infty $,
there exists a solution $F$ to  the Poisson equation $F-PF=f
-\langle  \pi, f \rangle$ such that $\langle
\pi, F^2 \rangle<\infty $. 
We consider the Boltzmann  case: the
selection kernel $\kappa$ is given by \reff{eq:Bgen}. For
 $\beta(x,A,y)$ respectively equal to $F(y)$ and $f(y)$, one has 
\[
\sigma(f,F)^2=\frac{1}{2}\Big(\sigma(f)^2-\Var_\pi (f) \Big)\mbox{
  and }\sigma(f,f)^2=\sigma(f)^2-\Delta(f).
\]
The
non-negative term $\Delta(f)$ is positive when $\Var_\pi(f)>0$.

Furthermore, if $\Var_\pi(f)>0$, then $\langle \pi,  f^2-fPf\rangle=\frac{1}{2}
\E_\pi \left[ (f(X_0)-f(X_1))^2\right]$ is positive,  the function
  $b\mapsto \sigma(f, bf)^2$ is minimal at 
\begin{equation}
   \label{eq:bstar2}
b_\star=\frac{\langle \pi, f^2 \rangle -\langle \pi, f \rangle^2}{\langle
  \pi, f^2 -fPf \rangle}, 
\end{equation}
and $b_\star>1$ when $\sigma(f,f)^2>0$.
\end{prop}

\begin{proof}
 Recall notations from Example \ref{exa:BMgen}.  We set
 $\kappa^B_A(dy)=\kappa^B(A,dy)$. For  $g$ and $h$  real valued 
  functions defined on $E$, we have
\begin{align}
\label{kkb}
\int \pi(dx)\cQ(x,dA)
\; \langle \kappa^B_A, g\rangle \langle \kappa^B_A,
h \rangle 
&=\int N_0(dA) r_A(dx) \; \langle \kappa^B_A, g\rangle \langle \kappa^B_A,
h \rangle \\ 
\nonumber
&=\int N_0(dA) \;  \langle r_A, g\rangle \langle \kappa^B_A,
h \rangle \\ 
\nonumber
&=\int \pi(dx) \cQ(x, dA) \; g(x)  \langle \kappa^B_A,
h \rangle \\ 
\nonumber
&=\langle \pi,gPh\rangle, 
\end{align}
where we used \reff{eq:Bgen} for the second equality. 
Using this equality with $h=g=F$ in the first term of the expression of
$\sigma(f,F)^2$ given in \reff{eq:sfF2}, we obtain
\[
\sigma(f,F)^2=\langle\pi,FPF-(PF)^2\rangle
=\frac{1}{2}\langle\pi,F^2-(PF)^2-(F-PF)^2\rangle
=\frac{1}{2}(\sigma(f)^2-\Var_\pi(f)),
\]
where we used the Poisson equation \reff{Poisson} for the last equality.

We also  get that
\begin{multline*}
   \int
 \pi(dx)\cQ(x,dA)\left[\Var_{\kappa^B_A} (bf-F) 
  - \Var_{\kappa^B_A} (F)\right]\\
\begin{aligned}
&=  \int
 \pi(dx)\cQ(x,dA)\left[\langle \kappa^B_A , (bf-F)^2\rangle -\langle
   \kappa^B_A , bf\rangle^2+2\langle
   \kappa^B_A , bf\rangle\langle
   \kappa^B_A , F\rangle \right.\\
&\hspace{6cm} \left. -\langle
   \kappa^B_A , F\rangle^2- \langle \kappa^B_A , F^2\rangle +
   \langle \kappa^B_A , F\rangle ^2\right]\\
&=  \langle \pi,  b^2f^2-2bfF-b^2fPf+2fPF\rangle\\
&=  b^2 \langle \pi,  f^2-fPf\rangle -2b \left(\langle \pi, f^2\rangle -
\langle \pi, f \rangle^2\right), 
\end{aligned}
\end{multline*}
where we used \reff{kkb} for  the second equation and \reff{Poisson} for
the last  equality. We  deduce from \reff{varFdel}  with $\beta_{x,A}=bf$
that
\[
\sigma(f,bf)^2 -\sigma(f)^2= b^2 \langle \pi,  f^2-fPf\rangle -2b
\left(\langle \pi, f^2\rangle - 
\langle \pi, f \rangle^2\right).
\]

We first check that $\Var_\pi(f)>0$ implies that $\langle \pi,
f^2-fPf\rangle >0$. 
If, when $X_0$ is distributed according to $\pi$, a.s. $f(X_1)=f(X_0)$,
then a.s. $k\mapsto f(X_k)$ is constant and by the ergodic theorem this
constant is equal to $\langle  \pi, f\rangle$. Therefore $\Var_\pi(f)>0$
implies positivity of $\langle  \pi, f^2-fPf\rangle$ which is equal to
$\frac{1}{2}\E_\pi  \left[ (f(X_0)-f(X_1))^2\right]$ by reversibility of
$\pi$ w.r.t. $P$.

Hence when $\Var_\pi(f)>0$, then  $b\mapsto \sigma(f,bf)^2$ is minimal
for $b=b_\star$ defined by \reff{eq:bstar2}. 

For the choice $b=1$, one obtains
\begin{equation}
   -\sigma(f,f)^2+\sigma(f)^2
= \langle \pi, f(f+Pf)\rangle-2\langle \pi, f \rangle^2 = \Delta(f)=\Var_\pi(f)+\langle \pi,
f_0Pf_0\rangle.\label{eqdifsffsf} 
\end{equation}
By \reff{kkb}, $\langle \pi, f_0Pf_0 \rangle=\int \pi(dx)
\cQ(x,dA)  \langle  \kappa_A^B,f_0 \rangle^2\geq 0$ and $\Delta(f)$ is
positive when $\Var_\pi(f)>0$. 

Moreover the difference $\langle \pi, fPf
 \rangle - \langle \pi, f \rangle^2=\langle \pi, f_0Pf_0 \rangle$ is
 non-negative thanks to \reff{kkb}
and when it is equal to $0$, then \reff{kkb} implies that $\langle
 \pi, f_0Pg \rangle=\langle \pi, gPf_0 \rangle=0$ for each function $g$
 on $E$ such that $\langle \pi,g^2\rangle<+\infty$. In this case, by
 \eqref{eqdifsffsf}, 
\begin{align*}
   \sigma(f,f)^2&=\sigma(f)^2+\sigma(f,f)^2-\sigma(f)^2\\
&=\langle
   \pi,(F+PF)(F-PF)\rangle-\Var_\pi(f)\\&=\langle
   \pi,(f_0+2PF)f_0\rangle-\Var_\pi(f)=0.\end{align*}
Hence when $\Var_\pi(f)>0$ and $\sigma(f,f)^2>0$ then, we have  $\langle \pi,
f_0Pf_0 \rangle>0$ and $b_\star>1$. 
\end{proof}
\section{Further results in the single-proposal case}\label{sp}
The Metropolis-Hastings  algorithm corresponds to the single proposal
case that is  the particular case of
the multi-proposal algorithm of Section \ref{sec:Egene} where  $\cQ(x,.)$
gives full weight to the set of subsets of $E$ (not assumed to be finite) containing $x$ and at
most one other element of $E$. The acceptance probability is then given
by $\rho(x,y)=\kappa(x,\{x,y\},y)$ and the selection kernel $Q(x,.)$ is
the image of $\cQ(x,.)$ 
by any measurable mapping such that the image of $\{x,y\}$ is $y$.
See Remark \reff{rem:multi-1} in the particular case of $E$ finite.
Equation \eqref{eq:reversmult2} is then equivalent to the following
generalization of \eqref{prerevers}
\begin{equation}
  \pi(dx)Q(x,dy)\rho(x,y)=\pi(dy)Q(y,dx)\rho(y,x).
\label{prereversgen}\end{equation}
Moreover the transition kernel of the Markov chain $X$ is given by
\begin{equation}
   \ind_{\{y\neq x\}}P(x,dy)=\ind_{\{y\neq x\}}\rho(x,y)Q(x,dy)\mbox{ and
}P(x,\{x\})=1-\int_{z\neq x}\rho(x,z)Q(x,dz).\label{noysinggen}
\end{equation}

Motivated by the study of the WR algorithm which corresponds to $\psi=f$
and of  the optimal choice $\psi=F$,  we are first going  to derive more
convenient  expressions of  $\sigma(f,\psi)^2$  in the
single  proposal framework.   We then  use this  new expression   to
construct a counter-example such that $\sigma(f,f)^2>\sigma(f)^2$. 
And, when 
$\rho(x,y)+\rho(y,x)$ is constant on $E^2_*=E^2\setminus\{(x,x):x\in
E\}$, using again the expression  of  $\sigma(f,\psi)^2$, we  compute
the  value of $b$ such that $\sigma(f,bf)^2$ is minimal and check that 
$\sigma(f,f)^2<\sigma(f)^2$ as soon as $f$ is non constant. 

\subsection{Another expression of the asympotic variance}
We recall that in the  notation $\E_\pi$, the subscript $\pi$ means that
$X_0$ is distributed according to $\pi$.
\begin{lem}
   \label{lem:added}
We assume that  $\langle \pi, f^2 \rangle<\infty $ and 
there exists a solution $F$ to  the Poisson  equation
  \eqref{Poisson}  such  that  $\langle\pi,F^2\rangle<+\infty$. Let
  $\psi$ be  square integrable:  $\langle \pi, \psi^2 \rangle<\infty $. 
 In  the
  single proposal case, we have
\begin{multline*}
\sigma(f,\psi)=\sigma(f)^2 -
 \E_\pi\left[\Big(1-\rho(X_0,
  X_1)\Big)\Big(F(X_1)-F(X_0)\Big)^2\right]\\
+ \E_\pi\left[\Big(1-\rho(X_0,
  X_1)\Big)\Big(\psi(X_1) -F(X_1)-\psi(X_0)+F(X_0)\Big)^2\right].
\end{multline*}
\end{lem}

\begin{proof}
In the single proposal case, $\kappa(x,\{x,y\},y)= 1-
  \kappa(x,\{x,y\},x)=\rho(x,y)$ for $x\neq y$. Therefore, for a real valued function $g$
  defined on $E$, we have
\begin{align}
\label{covber} 
\Var_{\kappa(x,\{x,y\},.)}(g)=\rho(x,y)(1-\rho(x,y))(g(y)-g(x))^2.
\end{align}
Thus we deduce that
\begin{align*}
 \int_{E^2_* }\pi(dx)Q(x,dy)\Var_{\kappa(x,\{x,y\},.)}(g)
=& \int_{E^2_* }\pi(dx)Q(x,dy)\rho(x,y)(1-\rho(x,y))
(g(y)-g(x))^2\\
=& \int_{E^2_* }\pi(dx)P(x,dy)(1-\rho(x,y))(g(y)-g(x))^2\\
=&\E_\pi\left[\Big(1-\rho(X_0,X_1)\Big)\Big(g(X_1)-g(X_0)\Big)^2\right].
\end{align*}
where we  used \eqref{noysinggen} for  the second  equality.
Plugging this formula with $g=\psi-F$ and $g=F$ in \eqref{varFdel} gives
the result.    
\end{proof}
Taking $\psi=F$ and $\psi=f$ in the previous Lemma gives the following
Corollary. 
\begin{cor}\label{singprop} 
We assume that  $\langle \pi, f^2 \rangle<\infty $ and there  exists  a
solution $F$  to  the Poisson  equation 
  \eqref{Poisson}  such  that  $\langle\pi,F^2\rangle<+\infty$.  In  the
  single proposal case, we have:
\[
\sigma(f,F)^2-\sigma(f)^2
=-\E_\pi\left[(1-\rho(X_0,X_1))(F(X_1)-F(X_0))^2\right],
\]
\[
   \sigma(f,f)^2-\sigma(f)^2 
=-\E_\pi\left[(1-\rho(X_0,X_1))
  \left[(F(X_1)-F(X_0))^2-(PF(X_1)-PF(X_0))^2\right]\right] . 
\]
\end{cor}

\subsection{A counter-example}
\label{sec:ce}
We are going to construct a counter-example such that
$\sigma(f,f)^2>\sigma(f)^2$ in the Metropolis case, thus
proving the statements concerning this case in Proposition \ref{metbark}. This counter-example is also such that the optimal
choice $\psi=F$ does not achieve variance reduction : $\sigma(f,F)^2=\sigma(f)^2$. 
Let $P$ be an irreducible transition matrix on $E=\{a,b,c\}$, with
invariant probability measure $\pi$ s.t. $P$ is reversible w.r.t. $\pi$, 
$$P(a,b)>0,\;P(a,a)>0\mbox{ and }P(a,c)\neq P(b,c).$$
Let $f$ be defined by $f(x)=\ind_{\{x=c\}} -P(x,c)$ for $x\in E$. We have
$$\langle \pi, f \rangle =\pi(c) -\sum_{x\in E}\pi(x) P(x,c)=0.$$ The function
$F(x)=\ind_{\{ x=c\}}$  solves the Poisson equation \reff{Poisson}:
$F-PF=f-\langle \pi, f \rangle$. 

Let  $\rho\in\left(\frac{P(a,b)}{P(a,a)+P(a,b)},1\right)$. We set
\begin{equation*}
   Q(x,y)=\begin{cases}\frac{P(a,b)}{\rho}\mbox{ if }(x,y)=(a,b),\\
P(a,a)-P(a,b)(\frac{1}{\rho}-1)\mbox{ if }(x,y)=(a,a),\\
P(x,y)\mbox{ otherwise.}
\end{cases}
\end{equation*}
We choose 
\begin{equation*}
   \rho(x,y)=\begin{cases}\rho\mbox{ if }(x,y)=(a,b),\\
1\mbox{ otherwise.}
\end{cases}
\end{equation*}
Since $\rho(a,b)\pi(a)Q(a,b)=\rho\pi(a)P(a,b)/\rho$, we have
$\rho(x,y)\pi(x)Q(x,y)=\pi(x)P(x,y)$ for all $x\neq y\in E$. Equation
\eqref{prerevers} follows from the reversibility of $\pi$ for $P$.
Notice also that \reff{casga} holds with $\gamma(u)=\min(1,u)$.

By construction, the matrix $P$ satisfies \reff{eq:def-P}. By
Corollary \ref{singprop}, we have $\sigma(f,F)^2-\sigma(f)^2=0$ and
\begin{equation}
   \label{eq:sff>sf}
\sigma(f,f)^2-\sigma(f)^2=\pi(a)P(a,b)(1-\rho)(P(b,c)-P(a,c))^2>0.
\end{equation}


Let us illustrate these results by simulation for the following specific
choice 
$$\pi=\frac{1}{10}\left(\begin{array}{c}6 \\ 3 \\1\end{array}\right),\;P=\frac{1}{60}\left(\begin{array}{ccc}38 & 21 & 1\\ 42 & 0 & 18\\ 6 &
    54 &
    0\end{array}\right),\;\rho=\frac{4}{10}\mbox{ and }Q=\frac{1}{120}\left(\begin{array}{ccc}13 & 105& 2\\84 & 0 &36\\ 12 & 108 & 0
   \\
\end{array}\right).$$
Then
$\sigma(f)^2-\sigma(f,f)^2=-0.010115$ amounts to $14\%$ of
$\sigma(f)^2\simeq 0.0728333$.

Using  $N=10\;000$ simulations, we give estimations of the variances
$\sigma_n^2$ of $I_n(f)$,  $\sigma^2_{WR,n}$ of $I_n(f,f)$ and of the
difference $\sigma_n^2-\sigma^2_{WR,n}$ with
asymptotic confidence intervals at level 95\%. The initial variable $X_0$
is generated according to the reversible probability measure $\pi$.

\begin{center}
   \begin{tabular}{|c|c|c|c|}
\hline
$n$ &
$\sigma_n^2$&
$\sigma_{WR,n}^2$&
$\sigma_n^2- \sigma_{WR,n}^2$\\
\hline
1 &  [0.1213 ,  0.1339]&    [0.1116 ,  0.1241]&[0.0091 ,  0.0104]\\
\hline
2 & [0.0728 ,  0.0779]&    [0.0758 ,  0.0815]&[-0.0041 , -0.0025]\\
\hline
5 &  [0.0733 ,  0.0791]&   [0.0798 ,  0.0859]&[-0.0075 , -0.0058]\\
\hline
10 &   [0.0718 ,  0.0772]& [0.0800 ,  0.0859]&[-0.0094 , -0.0074]\\
\hline
100 &     [0.0702 ,  0.0751]&   [0.0803 ,  0.0858]&[-0.0114 , -0.0092]\\
\hline
1000 &   [0.0719 ,  0.0769] &   [0.0811 ,  0.0867]&[-0.0105 , -0.0083]\\
\hline
\end{tabular}

\end{center}
\subsection{Case of a constant sum $\rho(x,y)+\rho(y,x)$.}

Under Boltzmann  selection rule, according to Proposition \ref{propbarker}, the  asymptotic variance $\sigma(f,f)^2$  of $I^{WR}_n(f)=I_n(f,f)$
is  smaller  than   the  one  $\sigma(f)^2$  of  $I_n(f)$ and
$\sigma(f,bf)$ is minimal for $b=b_\star$ given by \reff{eq:bstar2}. 
In   the  single   proposal  case,  Boltzmann  selection rule ensures  that
$\rho(x,y)+\rho(y,x)=1$ on $E^2_*=E^2\setminus\{(x,x):x\in E\}$. It
turns out that  we are still able  to prove the same results as  soon  as
$\rho(x,y)+\rho(y,x)$ is  constant on $E^2_*$. Notice that 
$\Var_\pi(f)\geq 0$ and that the trivial case  $\Var_\pi(f)=0$
corresponds to $f$ constant $\pi$-a.s.. 
\begin{prop}
  \label{prop:geneB} We assume $\langle \pi, f^2 \rangle<\infty $,
  $\Var_\pi(f)>0$, there  exists a solution $F$ to  the Poisson equation
  \eqref{Poisson}  such that  $\langle  \pi, F^2  \rangle<\infty $.   We
  consider the single proposal case  and assume that there exists $\alpha \in
  (0,2)$ such that
\begin{equation}
\pi(dx)Q(x,dy)\mbox{ a.e. on }E^2_*,\;\rho(x,y)+\rho(y,x)=\alpha. 
\label{consts}\end{equation}
Then we have:
\begin{itemize}
   \item[i)] $\langle \pi,  f^2-fPf\rangle=
\frac{1}{2}\E_\pi \left[ (f(X_0)-f(X_1))^2\right]$ is positive.
\item[ii)]\begin{align}
\label{eq:alpha-psi}
\sigma(f,\psi)-\sigma(f)^2 =&-(1-\alpha/2)
 \E_\pi\left[\Big(F(X_1)-F(X_0)\Big)^2\right]\\
&+(1-\alpha/2)  \E_\pi\left[\Big(\psi(X_1)
  -F(X_1)-\psi(X_0)+F(X_0)\Big)^2\right]\notag,   
\end{align} for any real valued function $\psi$ on $E$ such
  that $\langle \pi, \psi^2 \rangle<\infty$.

   \item[iii)] The function
  $b\mapsto \sigma(f, bf)^2$ is minimal at $b_\star$ given by
  \reff{eq:bstar2} and $b_\star\geq 1/\alpha$. 
   \item [iv)] $\sigma(f,f)^2-\sigma(f)^2=-(2-\alpha )\Delta(f)<0$,
where $\Delta(f)$ is given by \reff{eq:Deltaf}. 
\end{itemize}  
\end{prop}

\begin{proof}
Statement i) follows from the proof of Proposition \ref{propbarker}. 

For statement ii), notice that by reversibility of $\pi$, we deduce from
Lemma \ref{lem:added} that 
\begin{align*}
\sigma(f,\psi)-\sigma(f)^2
   &=-
 \E_\pi\left[\Big(1-\rho(X_1,
  X_0)\Big)\Big(F(X_1)-F(X_0)\Big)^2\right]\\
&\hspace{1cm}+ \E_\pi\left[\Big(1-\rho(X_1,
  X_0)\Big)\Big(\psi(X_1) -F(X_1)-\psi(X_0)+F(X_0)\Big)^2\right]. 
\end{align*}
This and Lemma \ref{lem:added} imply \reff{eq:alpha-psi}. 

For   iii),   using    \reff{eq:alpha-psi}   with   $\psi=bf$,   it   is
straightforward to get that $\sigma(f, bf)^2$ is minimal when $b$ equals
\[
\frac{\E_\pi\left[(f(X_1) -f(X_0)) (F(X_1)-F(X_0))\right]}{\E_\pi\left[(f(X_1)
 -f(X_0))^2\right]}=\frac{\langle \pi, f(F-PF) \rangle}{\langle \pi, f^2
-fPf\rangle}=\frac{\langle \pi, f^2 \rangle -\langle \pi, f \rangle^2}{\langle
  \pi, f^2 -fPf \rangle}=b_\star.
\]
Remarking that $ b_\star=\frac{\langle \pi, f_0^2 \rangle}{\langle \pi,
  f_0^2-f_0Pf_0\rangle}=\frac{\langle \pi, f_0^2 \rangle}{\alpha\langle
  \pi, f_0^2\rangle-\langle\pi,f_0Pf_0+(\alpha-1)f_0^2\rangle}$ and
using Lemma \ref{lem:hPh} below,  one
deduce  that $b_\star\geq 1/\alpha$.

We now prove  iv). Recall that  $f_0=f-\langle \pi,f  \rangle$. 
Since $\langle \pi, f_0(f_0+ Pf_0) \rangle=(2-\alpha)\Var_\pi(f)+\langle
\pi, f_0Pf_0+(\alpha-1)f_0^2\rangle$, 
we deduce from Lemma \ref{lem:hPh} that $\Delta(f)$ given by
\reff{eq:Deltaf}
is positive. We have
\begin{multline*}
   \inv{2} \E_\pi\left[(f(X_1)-F(X_1)-f(X_0)+F(X_0))^2
     -(F(X_1)-F(X_0))^2\right]\\ 
\begin{aligned}
&=   \inv{2} \E_\pi\left[(f_0(X_1)-f_0(X_0))^2 \right] -
\E_\pi\left[(f_0(X_1)-f_0(X_0))( F(X_1) -F(X_0)) \right]\\
&=\langle \pi, f_0^2-f_0Pf_0\rangle -2\langle \pi, f_0(F-PF) \rangle\\
&=-\langle \pi, f_0(f_0+ Pf_0) \rangle,
\end{aligned}
\end{multline*}
where we used that  $\pi$ is invariant for $P$ and that $P$ is reversible
with respect to $\pi$ for the second equality and that $F$ solves
\eqref{Poisson} for the last equality. We conclude using  
\reff{eq:alpha-psi} with $\psi=f$.
\end{proof}

\begin{lem}
   \label{lem:hPh}
   Let $h$ be a real valued  function defined on $E$ such that $\langle
   \pi, h^2 \rangle<\infty $.  
Under hypothesis \reff{consts}, we have $\langle \pi, hPh +(\alpha-1)  h^2
\rangle\geq 0$. 
\end{lem}
\begin{proof}

Using \eqref{noysinggen} then
\eqref{consts}, we obtain
\begin{align*}
 \langle \pi, hPh +(\alpha-1) h^2 \rangle
&=\int_{E^2_*}\pi(dx)Q(x,dy)\rho(x,y)h(x)h(y)\\
&\hspace{2cm}+\int_E\pi(dx)\left(\alpha-\int_{E}\ind_{y\neq x}Q(x,dy)\rho(x,y)\right)h^2(x)\\
&=\int_{E^2_*}\pi(dx)Q(x,dy)\left[\rho(x,y)h(x)h(y)+\rho(y,x)h^2(x)\right]\\
&\hspace{2cm}+\alpha\int_E\pi(dx)Q(x,\{x\})h^2(x).
\end{align*}
To conclude, it is enough to check that the first term in the r.h.s. is
nonnegative. Using \eqref{consts} and  \eqref{prereversgen} for the
first equality, we get 
\begin{align*}
\alpha \int_{E^2_*}\pi(dx)Q(x,dy)
&\left[\rho(x,y)h(x)h(y)+\rho(y,x)h^2(x)\right]\\
&=\int_{E^2_*}\pi(dx)Q(x,dy)\rho(y,x)\left[\rho(x,y)h(x)h(y)
  +\rho(y,x)h^2(x)\right]\\ 
&\hspace{2cm}+\int_{E^2_*}\pi(dy)Q(y,dx)\rho(y,x)
\left[\rho(x,y)h(x)h(y)+\rho(y,x)h^2(x)\right]\\
&=\int_{E^2_*}\pi(dx)Q(x,dy)\left[\rho(y,x)h(x)+\rho(x,y)h(y)\right]^2\\
&\geq 0.
\end{align*}
   
\end{proof}

\section{Other remarks}
We  work in  the general  setting  of Section  \ref{sec:Egene}.

\subsection{About the estimator $I_n(f+ P\psi - \psi)$}
\label{sec:compvarppp}
Motivated by Remark \ref{rem:Reduc} on the study of $I_n(f+P\psi
-\psi)$, we compute the asymptotic variance $\tilde\sigma(f,\beta)^2$ of
\begin{multline*}
   I_n(f)+\inv{n}
   \sum_{k=0}^{n-1}\;\;\left(\int 
     \cQ(X_k,dA)\kappa(X_k,A,d\tilde{x}) \beta(X_k, A, \tilde x)  - 
 \beta(X_k,A_{k+1},X_{k+1})\right)\\
=I_n(f) + \inv{n} 
   \sum_{k=0}^{n-1} \bigg( \E[\beta(X_k, A_{k+1}, X_{k+1})|X_k] -
   \beta(X_k,A_{k+1},X_{k+1})\bigg). 
\end{multline*}

Following the proof of Theorem \ref{theo:princ2}, one obtains that the above
estimator of $\langle \pi, f \rangle$ is under the hypotheses of Theorem
\ref{theo:princ2} convergent and asymptotically normal with asymptotic
variance 
\[
\tilde\sigma(f,\beta)^2=\sigma(f,\beta)^2+
\int\pi(dx)\left[\Var_{\cQ(x,.)}(\kappa\beta_x
  -\kappa F_x)-\Var_{\cQ(x,.)}(\kappa F_x)\right], 
\] 
where $\displaystyle \Var_{\cQ(x,\cdot)} (\varphi)=\int \cq(x, dA)
\varphi(A)^2 - 
\left(\int \cq(x, dA) \varphi(A)\right)^2 $, $\kappa
\beta_x(A)=\langle \kappa_{x,A} ,\beta_{x,A} \rangle$ and $\kappa
F_x(A)=\langle \kappa_{x,A} ,F \rangle$. 

Notice that the sign of  $\tilde\sigma(f,\beta)^2-\sigma(f,\beta)^2$ depends
on $\beta$ (take $\beta_{x,A}=F$ and $\beta_{x,A}=-F$).

\subsection{Changing the selection kernel in $\cJ_n$}
\label{sec:changingK}
Let $\kappa'\neq  \kappa$ be such that  \reff{eq:reversmult2} (or simply
\eqref{reversmult}  if  $E$ is  finite)  still  holds  when $\kappa$  is
replaced by $\kappa'$ and  $\cJ'_n(\psi)$ and $\cJ_n'(\beta)$ be defined
like $\cJ_n(\psi)$  and $\cJ_n(\beta)$ with the chain  $X$ unchanged but
with         $\kappa(X_k,A_{k+1},\tilde{x})$         replaced         by
$\kappa'(X_k,A_{k+1},\tilde{x})$ in \eqref{cjpsi} and \eqref{cjbet}.
Thus, we have 
\[
  \cJ_n'(\psi)=\inv{n} \sum_{k=0}^{n-1}\;\; \sum_{\tilde{x}\in
  A_{k+1}}\left(\kappa'(X_k,A_{k+1},\tilde{x})-\ind_{\{X_{k+1}
    =\tilde{x}\}}\right)\psi(\tilde{x})  .
\] 
Note that in general $ \sum_{\tilde{x}\in
  A_{k+1}}\kappa'(X_k,A_{k+1},\tilde{x})\psi(\tilde{x}) \neq
  \E[\psi(X_{k+1})|X_k, A_{k+1}]$. 

In  the  single  proposal  case,  Frenkel  \cite{f:wrmc}  suggests  that
$\cJ'_n(f)$ can  also be used as  a control variate.  In  general, for a
real  valued function  $\beta$ defined  on $E\times  \cp \times  E$, the
almost sure limit of $\cJ'_n(\beta)$ is different from zero, which means
the estimator $I_n  (f) + \cJ_n'(\beta)$ of $\langle  \pi, f \rangle$ is
not  convergent.  However,  when $\beta(x,A,  \cdot)=\psi(\cdot)$, Lemma
\ref{lem:Ip} below  ensures that the estimator $I_n  (f) + \cJ_n'(\psi)$
of $\langle  \pi, f \rangle$  is convergent.  It  is also easy  to prove
that this estimator is  asymptotically normal and compute the asymptotic
variance, but  we have not been  able to compare it  with the asymptotic
variance $\sigma(f)^2$ of $I_n(f)$.

\begin{lem}
\label{lem:Ip}
  We assume $X$ is Harris recurrent, $\langle \pi, f^2 \rangle<\infty $,
there exists a solution $F$ to  the Poisson equation $F-PF=f
-\langle  \pi, f \rangle$ such that $\langle
\pi, F^2 \rangle<\infty $,  and $\psi$ is such that: $\langle \pi,
\psi^2  \rangle <\infty $. Under those  assumptions, the estimator
$I_n(f)+ \cJ_n'(\psi) $ of $\langle \pi, f \rangle$  is  consistent:
a.s. $\displaystyle  \lim _{n\rightarrow \infty } I_n(f)+ \cJ_n'(\psi)=\langle
\pi, f \rangle$.   
\end{lem}

\begin{proof}
   We set
\[
\Delta R_n=  \int \kappa'(X_{n-1},A_n,d\tilde{x})\psi(\tilde x)- \int
\cQ(X_{n-1}, dA) 
\kappa'(X_{n-1},A,d\tilde{x})\psi(\tilde x).
\] 
Notice that $\Delta R_n$ is  square integrable and  that $\displaystyle
\E[\Delta  R_{n+1}|\cg_n]=0$,  where  $  \cg_n$  is  the  $\sigma$-field
generated by $X_0$ and $(A_i,  X_i)$ for $1\leq i\leq n$.  In particular
$R=(R_n, n\geq  0)$ with $R_n=\sum_{k=1}^n  \Delta R_k$ is  a martingale
w.r.t.  to the filtration $(\cg_n, n\geq 0)$. Notice that 
\[
\cJ_n'(\psi)= \inv{n} R_n + I_n(\gamma)-\inv{n}\int
\cQ(X_{n}, dA) 
\kappa'(X_{n},A,d\tilde{x})\psi(\tilde x) +\inv{n}\int
\cQ(X_{0}, dA) 
\kappa'(X_{0},A,d\tilde{x})\psi(\tilde x) ,
\]
where $\displaystyle \gamma(x)=\int \cQ(x,dA)
\kappa'(x,A,d\tilde x) \psi(\tilde x)-\psi(x) $. 
Following the proof of Theorem \ref{theo:princ2}, we easily get that
a.s.  $\displaystyle
\lim_{n\rightarrow  \infty } 
\inv{n} R _n =0$ and  that a.s. 
 \[
 \lim_{n\rightarrow\infty  }  \cJ_n'(\psi)  =\lim_{n\rightarrow\infty  }
 I_n(\gamma)=\langle \pi, \gamma \rangle.
\]
Using \reff{eq:reversmult2} satisfied by $\kappa'$ instead of $\kappa$,
we get that $\langle \pi, \gamma \rangle=0$. This ends the proof of the
Lemma. 
\end{proof}

\newcommand{\sortnoop}[1]{}


\begin{thebibliography}{1}

\bibitem{ad:hrud}
H. C. Andersen and P. Diaconis.
\newblock{Hit and Run as  a unifying device}
\newblock{\em J. de la SFdS et revue de stat. appli.} 148(4):5-28,2007. 

\bibitem{atchperr}
Y.F. Atchad\'e and F. Perron.
\newblock{Improving on the Independent Metropolis-Hastings algorithm}.
\newblock{\em Statist. Sinica} 15(1):3-18, 2005

\bibitem{a}
M.~Ath\`enes.
\newblock {Web ensemble averages for retrieving relevant information
  from rejected Monte Carlo moves}.
\newblock {\em Europ. Phys. J. B}, 58:83-95 (2007). .

\bibitem{cck}
D.~Ceperley, G.V. Chester and M.H. Kalos.
\newblock {Monte Carlo simulation of a many fermion study}.
\newblock {\em Phys. Rev. B} 16(7):3081-3099, 1977.

\bibitem{d:rim}
M.~Duflo.
\newblock {\em Random iterative models}, volume~34 of {\em Applications of
  Mathematics (New York)}.
\newblock Springer-Verlag, Berlin, 1997.
\newblock Translated from the 1990 French original by Stephen S. Wilson and
  revised by the author.

\bibitem{f1:wrmc}
D.~Frenkel.
\newblock {Speed-up of Monte Carlo simulations by sampling of rejected states}.
\newblock {\em Proc. Nat. Acad. Scienc.}, 101(51):17571-17575, 2004.

\bibitem{f:wrmc}
D.~Frenkel.
\newblock Waste-{R}ecycling {M}onte {C}arlo.
\newblock {\em  Lect. Notes in
  Phys.}, 703:127-137, Springer, 2006.

\bibitem{mt:mcss}
S.~P. Meyn and R.~L. Tweedie.
\newblock {\em Markov chains and stochastic stability}.
\newblock Communications and Control Engineering Series. Springer-Verlag London
  Ltd., London, 1993.

\bibitem{munos}
R. Munos.
\newblock {Geometric Variance Reduction in Markov Chains: Application to
Value Function and Gradient Estimation}.
\newblock {\em J. Machine Learning Res.},7:413-427, 2006.

\bibitem{p:omcsumc}
P.~H. Peskun.
\newblock Optimum {M}onte-{C}arlo sampling using {M}arkov chains.
\newblock {\em Biometrika}, 60:607--612, 1973.

\bibitem{cr:mcsm}
C.~P. Robert and G.~Casella.
\newblock {\em Monte {C}arlo statistical methods}.
\newblock Springer Texts in Statistics. Springer-Verlag, New York, 1999.


\end{thebibliography}
\end{document}